\documentclass[12pt]{amsart}
\usepackage{geometry} 
\usepackage[applemac]{inputenc}
\usepackage{textcomp}
\usepackage{amssymb}
\usepackage{graphicx}
\usepackage{subfigure}

\def\conv{\,\mbox{conv}\,}

\newtheorem{teorema}{Theorem}
\newtheorem{definicion}{Definition}

\newtheorem{problem}{Problem}
\newtheorem{lema}{Lemma}
\newtheorem{coro} {Corollary}

\title{Order types of convex bodies}
\author{
       Alfredo Hubard,
        Luis Montejano, Emiliano Mora, Andrew Suk
}

\date{} 

\begin{document}

\maketitle
\begin{center}\textbf{Abstract}\end{center}
We prove a \emph{Hadwiger transversal} type result, characterizing convex position on a family of non-crossing convex bodies in the plane. This theorem suggests a definition for the \emph{order type} of a family of convex bodies, generalizing the usual definition of order type for point sets. This order type turns out to be an oriented matroid. We also give new upper bounds on the Erd\H{o}s-Szekeres theorem in the context of convex bodies.
\section{Introduction}
In 1933 Paul Erd\H{o}s and  George Szekeres proved that \emph{for every $n \in \mathbb{N}$, there exists  $N(n) \in \mathbb{N}$ such that any set of at least $N(n)$ points in general position contains a subset of $n$ points forming a convex polygon} \cite{ErSze1}. They came out with two proofs, one attributed to Erd\H{o}s and one to Szekeres.
In 1989 Tibor Bisztriczky and G\'{a}bor Fejes T\'{o}th  generalized the classical Erd\H{o}s-Szekeres theorem to disjoint convex compact sets in the plane.
\begin{definicion} A family $\mathcal{F}$ of  sets in the plane    is \textbf{in convex position}  if \begin{displaymath} \textrm{ for every } A \in \mathcal{F} \quad  \conv(\mathcal{F}) \neq \conv(\mathcal{F}\setminus A) \end{displaymath}
\end{definicion}
\begin{definicion} A family of sets in the plane $\mathcal{F}$  is \textbf{in general position}  if every triple is in convex position.
\end{definicion}
The difficulty generalizing Szekeres's technique is that in contrast to the case of points there are families of convex bodies such that every four-tuple is in convex position, but the whole family is not. In fact there are families of $n$ convex bodies such that any $(n-1)$-tuple is in convex position, but the family fails to be in convex position.  The fundamental result of \cite{bis1} is the next theorem.

\begin{teorema}\cite{bis1}For every $n \in \mathbb{N}$, there is an $M(n)\in \mathbb{N}$ such that every disjoint  family $\mathcal{F}$  with at least $M(n)$  convex bodies in general position in the plane, contains a subfamily of $n$ sets in convex position.
\end{teorema}

Bisztriczky and G. Fejes T\'{o}th also made an appealing conjecture, that $N(n)=M(n)$.  In 1998 J\'{a}nos Pach and G\'{e}za T\'{o}th in \cite{pach1} improved their triply exponential upper bound on $M(n)$ by showing that
 $$M(n) \leq {2n-4 \choose n-2}^2+1.$$
Pach and T\'{o}th also showed in \cite{pach2} that if one relaxes the disjointness hypothesis to \emph{noncrossing}, the proposition still holds.

\begin{definicion} A family  $\mathcal{F}$ of sets in the plane is \textbf{noncrossing} if for every pair of convex bodies $\{A,B\}$, the sets $A \setminus B$ and $B\setminus A$ are connected\end{definicion}

\begin{teorema}\cite{pach2}For every $n \in \mathbb{N}$ there is a $M_{noncr}(n)\in \mathbb{N}$ such that every noncrossing family $\mathcal{F}$  with at least $M_{noncr}(n)$  convex compact sets in general position in the plane contains a subfamily of $n$ sets in convex position.
\end{teorema}

In this paper we give new proofs of these results yielding better bounds for $M_{noncr}(n)$ and $M(n)$. These proofs are applications of the next theorem which is our main result.

\begin{definicion}\label{orient} An ordered family $\mathcal{F}= \{A_1,A_2,A_3\}$ of three noncrossing convex bodies in the plane is said to be \textbf{clockwise oriented} (\textbf{counterclockwise oriented}) if there exist representative points $a_1 \in A_1 \cap bd[\conv(\mathcal{F})]$, $a_2 \in A_2 \cap bd[\conv(\mathcal{F})]$  and $a_3 \in A_3 \cap bd[\conv(\mathcal{F})]$ such that $a_1,a_2,a_3$ are oriented clockwise (counterclockwise).
\end{definicion}

\begin{teorema}\label{main} A family of noncrossing convex bodies is in convex position if and only if there exists an ordering of the family such that every triple is oriented counterclockwise.
 \end{teorema}

The last section of this paper discusses the connections between this transversal-Hadwiger type result and the Bisztriczky-Fejes T\'{o}th conjecture. We ask when a family of convex bodies is representable by points. To properly pose this problem we need to define \emph{order type} of a family of noncrossing convex bodies. It turns out that under natural assumptions the order type of a family of convex bodies is an oriented matroid. This new connection between two classical objects in discrete geometry suggests some new directions of research. We close this paper outlining these by posing some open problems.

\section{A new generalization}

Before going into the proof of Theorem \ref{main} we will prove a new generalization of the Erd\H{o}s-Szekeres Theorem for convex bodies. This generalization has the advantage over the previous ones of asking no conditions on the combinatorial geometry of the bodies, such as general position. It also suggests thinking of convex position as a \emph{transversal} property.  We denote the set of all $k$-tuples of  $X$ by ${X \choose k}$ and $R^3(k,l)$ will denote the Ramsey function for complete $3$-uniform hypergraphs, i.e. the minimal number such that if $|X|=R^3(k,l)$, then in every blue-red coloring of ${X \choose 3}$ there is a $Y \subset X$,  such that either $|Y|=k$ and every hyperedge $h \in {Y \choose 3}$ is red, or $|Y|=l$  and every hyperedge $h \in {Y \choose 3}$ is blue.

\begin{teorema}\label{convexo} For every pair $t,n \ge 3 $ there is an $h(t,n) \in \mathbb{N}$ such that  any planar family with more than $h(t,n)$ convex bodies contains either a subfamily of $t$ convex bodies with a common transversal line or a family of $n$ convex bodies in convex position.
 \end{teorema}

\noindent\textbf{Proof.} We will prove that $h(t,n) \le R^3(4t,N(n))$. First, color the triples and apply Ramsey's theorem. Let $X \in {\mathcal{F} \choose 3}$, color $X$ red if it has a transversal line, color $X$ blue if it does not have a transversal line. Ramsey's theorem yields a subfamily of $4t$ convex sets such that either each triple has a transversal line or a subfamily of $N(n)$ such that no triple has a transversal line. In the first case,  we may apply a result of J\"urgen Eckhoff

\begin{teorema} \cite{eckhoff} If $\mathcal{G}$ is a planar family of convex bodies such that every triple has a transversal line, then there are (at most) four lines such that every body in $\mathcal{G}$ is intersected by at least one of them.
\end{teorema}

  So we can conclude that there is a line that intersects at least $t$ of the bodies. In the second case, the subfamily that we obtained has no transversal line. Choosing a point in each set, we obtain $N(n)$ points in general position. By the Erd\H{o}s-Szekeres theorem there is a subset of $n$ points in convex position. The absence of transversal lines implies that the corresponding convex bodies will be in convex position. The last statement is not hard to prove directly; it also follows from Theorem \ref{main}. \begin{flushright}$ \blacksquare$\end{flushright}

\begin{coro}For every $t,n,d \ge 2 $, there is an $h_d(t,n) \in \mathbb{N}$ such that  any family with more than $h_d(t,n)$ convex bodies in $\mathbb{R}^d$ contains either $t$ members with a transversal hyperplane or $n$ members in convex position.
 \end{coro}

\noindent \textbf{Proof.}  We will prove $h_d(t,n) \le h(t,n)$. Apply Theorem \ref{convexo} to the image of any two-dimensional projection of the $d$ dimensional family. The pre-image of a planar family with a transversal line has a transversal hyperplane. The pre-image of a family in convex position is also in convex position. \begin{flushright}$ \blacksquare$\end{flushright}

With this corollary we obtain a family in convex position or a family with a transversal hyperplane; a natural open problem is to find for which values of $k$,  there is an $h_d^k(t,n)$ such that among $h_d^k(t,n)$ bodies in $\mathbb{R}^d$ there are either $n$ bodies in convex position, or $t$ bodies with a transversal $k$-flat. Note that for points, $h_d^2(t,n)$ exists for any (nontrivial) triple of natural numbers $d,t,n$.

 \section{A Hadwiger-type theorem}
Our previous Theorem suggests thinking of convex position as a transversal property.  Hadwiger's transversal theorem claims that a planar family of convex bodies has a transversal line if and only if there is an ordering of the family such that each triple has an oriented transversal line that intersects the sets in the prescribed ordering.\\
After the original proofs of \cite{ErSze1}, several other proofs of the Erd\H{o}s-Szekeres Theorem have been discovered. The next beautiful proof was posed as an exercise in \cite{grun}. Consider a point set in general position with more than $R^3(n,n)$ points in the plane. Order them in any way, and color a triple red if it is oriented clockwise and blue if it is oriented counterclockwise. By Ramsey's theorem there exists a subset of  $n$ points such that every triple is oriented likewise;  this implies that this $n$-set is in convex position.

The underlying geometric statement to prove convex position is the same that Erd\H{o}s used for the cups and caps technique, and is the same that Valtr and T\'{o}th used to obtain the best known upper bound on $N(n)$:  \emph{A set of points is in convex position if and only if there exists an ordering of the points such that every triple is oriented clockwise}.

\medskip

In the rest of the paper, we will assume that \emph{no two convex bodies are tangent, and no three convex bodies share a common tangent line}. We will work only in the plane and assume \emph{noncrossing} families. We restate the definition of orientation of such triples.

\vspace{.3cm}

\noindent\textbf{Definition \ref{orient}.} \emph{An ordered family $\mathcal{F}= \{A_1,A_2,A_3\}$ of three noncrossing convex bodies in the plane is said to be \textbf{counterclockwise oriented} if there exist representative points $a_1 \in A_1 \cap bd[\conv(\mathcal{F})]$, $a_2 \in A_2 \cap bd[\conv(\mathcal{F})]$  and $a_3 \in A_3 \cap bd[\conv(\mathcal{F})]$ such that $a_1,a_2,a_3$ are oriented counterclockwise.}

\vspace{.3cm}
\noindent \emph{Remark 1:} Note that an ordered triple of noncrossing convex bodies can have one, two, or no orientation.\\\\
\emph{Remark 2:} An ordered family of  noncrossing  convex bodies is in general position if and only if every triple has at least one orientation.\\\\
\emph{Remark 3:} An ordered triple has both orientations if, and only if, one of the bodies disconnects the convex hull of the triple. \\\\
We say that $X$ \emph{disconnects} $\mathcal{F}$ if $conv(\mathcal{F})\setminus X$ is disconnected.  Likewise we say that $\mathcal{F}$ is \emph{disconnectable} if there exists an $X \in \mathcal{F}$ such that $X$ disconnects $\mathcal{F}$.  See Figure 2(a).  Notice that a triple has both orientations if and only if it is disconnectable.  We are now ready to prove our main result.

\vspace{.3cm}

\noindent\textbf{Theorem \ref{main}.} \textit{A family of noncrossing convex bodies is in convex position if and only if there exists an ordering of the family such that every triple is oriented clockwise.}

\medskip

\noindent \textbf{Proof.} Let $\mathcal{F} = \{C_1,C_2,...,C_N\}$ be an ordered family of non-crossing convex sets in general position such that every triple has a clockwise orientation.  For sake of contradiction, suppose that there exists a convex body $C_j \in \mathcal{F}$ such that $C_j \subset conv(\mathcal{F}\setminus C_j)$.  Then let $\mathcal{F}_0\subset \mathcal{F}\setminus C_j$ be the minimum size subfamily such that $C_j \subset conv(\mathcal{F}_0)$.  By minimality we know that $\mathcal{F}_0$ is not disconnectable, and all members in $\mathcal{F}_0$ appear on the boundary of $conv(\mathcal{F}_0)$.  Let $\mathcal{F}_0 = \{C_{i_1},C_{i_2},...,C_{i_m}\}$ denote the order of the convex bodies as they appear in clockwise order along the boundary of $conv(\mathcal{F}_0)$.

\medskip

\noindent \textbf{Observation.}  For $1 \leq k < m$, \emph{$(C_{i_k},C_{i_{k + 1}},C_j)$ must have a unique orientation.  That is, $(C_{i_k},C_{i_{k + 1}},C_j)$ is not disconnectable.}

 \medskip

\noindent \textbf{Proof.}  Notice that $C_j$ cannot disconnect $conv(C_{i_k}\cup C_{i_{k + 1}}\cup C_j)$.  Assume that there exists a $k$ such that $C_{i_k}$ disconnects $conv(C_{i_k}\cup C_{i_{k + 1}} \cup C_j)$.  See Figure 2.b. Then this contradicts the minimality of $\mathcal{F}_0$ since this would imply $C_j \subset conv(\mathcal{F}_0\setminus C_{i_{k+1}})$.  Likewise, if $C_{i_{k + 1}}$ disconnects $conv(C_{i_k}\cup C_{i_{k + 1}} \cup C_j)$, then $C_j \subset conv(\mathcal{F}_0\setminus C_{i_{k}})$.

$\hfill\square$

 \begin{figure}
  \centering
    \subfigure[$X$ disconnects the family.]{\label{hide2}\includegraphics[width=.33\textwidth]{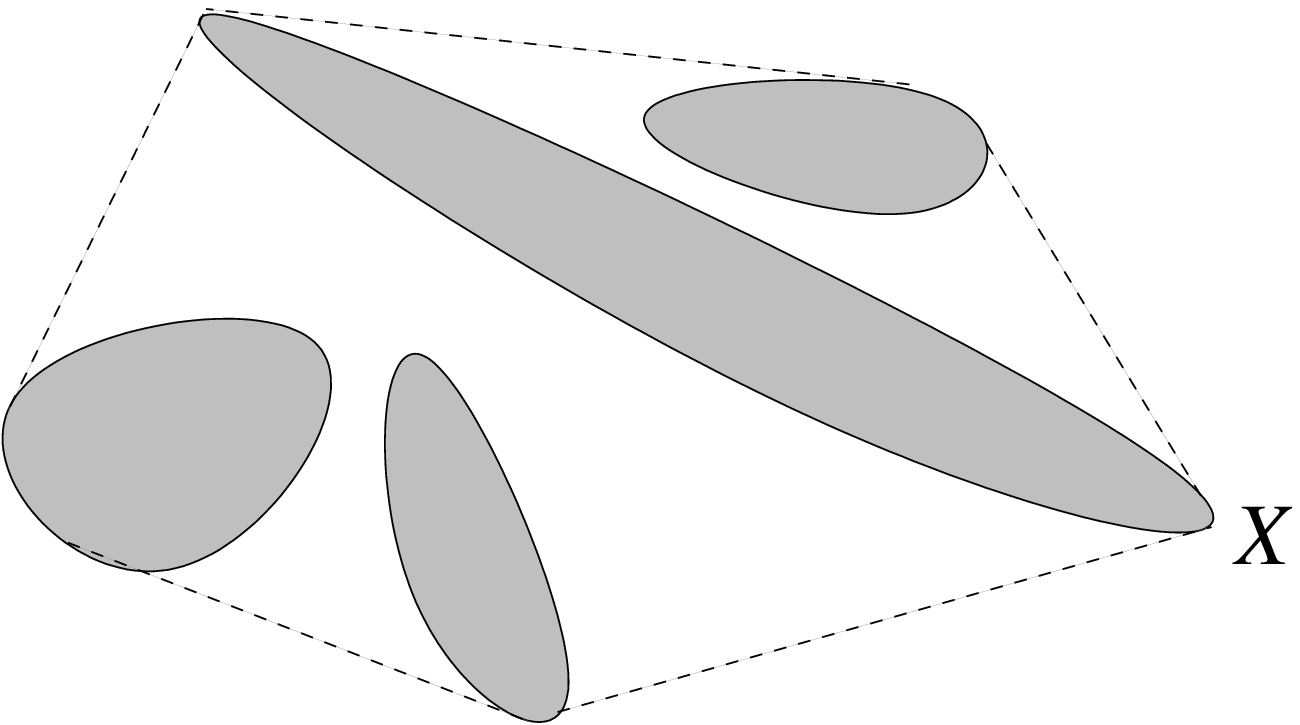}} \hspace{2cm}
      \subfigure[$C_{i_k}$ disconnects $conv(C_{i_k}\cup C_{i_{k + 1}} \cup C_j)$.]{\label{clockwise}\includegraphics[width=0.33\textwidth]{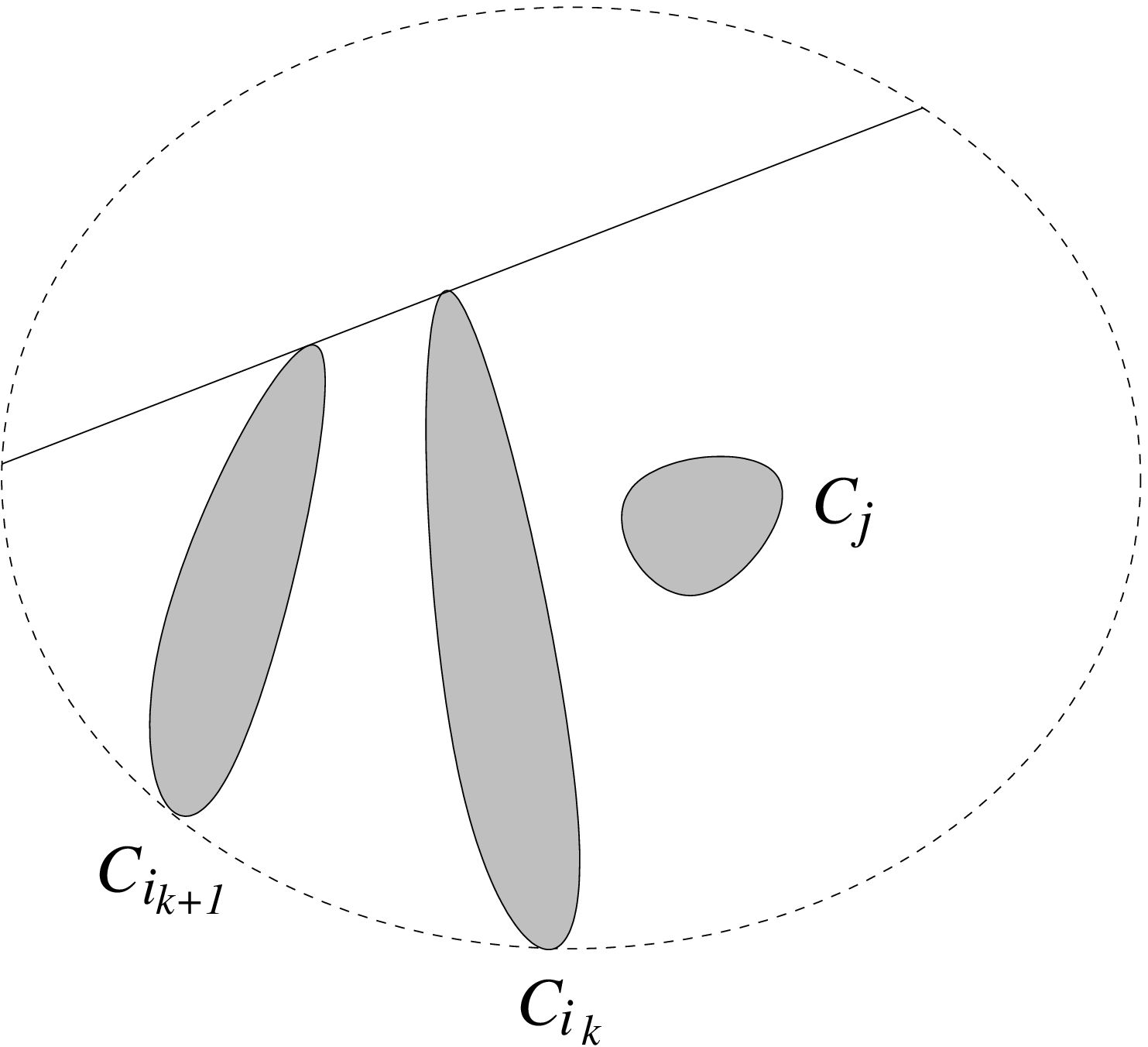}}

\caption{}
\label{fig:disconnect}
\end{figure}

\noindent Since $(C_{i_1},C_{i_2},C_j)$ only has a clockwise orientation, this implies we have the cyclic ordering $(i_1,i_2,j)$.  Likewise since $(C_{i_2},C_{i_3},C_j)$ only has a clockwise orientation, we have the cyclic ordering $(i_1,i_2,i_3,j)$.  As we continue around, we have the cyclic ordering $(i_1,i_2,i_3,...,i_k,j)$.  However this implies that $(C_{i_1},C_{i_{k}},C_j)$ has a counterclockwise orientation, and hence we have a contradiction.

\medskip

For the other direction suppose $\mathcal{F}$ is in convex position.  Then by starting at a point $p\in C_1$ that lies on the boundary of $conv(\mathcal{F})$, we order $\mathcal{F}$ as each body appears for the first time in clockwise order along the boundary of $conv(\mathcal{F})$.  Then every triple has a clockwise orientation.

$\hfill\square$

As an immediate Corollary, we have.

\begin{coro}$M_{noncr}(n)\le R^3(n,n)$. \end{coro}
The next Theorem improves the bound of \cite{pach1} by a factor of $2$.  The strategy that we use is a combination of the ones on \cite{pach1} and \cite{vato}. We will need the following hypergraph version of the \cite{ErSze1} and a Lemma by Pach and T\'oth.

 \begin{lema}\cite{ErSze1}Given a complete $3$-uniform hypergraph on ${k+l+4 \choose l+2}+1$  vertices. Assume there is an order on the vertices and $\chi$ a two coloring on the edges such that if $i<j<k<l$ and $\chi(x_i,x_j,x_k)=\chi(x_j,x_k,x_l)$ then $\chi(x_i,x_k,x_l)=\chi(x_i,x_j,x_l)=\chi(x_i,x_j,x_k)=\chi(x_j,x_k,x_l)$.  Then there is a complete blue subgraph on $l+1$ vertices or a complete red subgraph on $k+1$ vertices. \end{lema}

 \begin{lema} \cite{pach1} If $\mathcal{F}$ is a family of ${2n-4\choose n-2} + 1$ pairwise disjoint convex sets with a line transversal, then $\mathcal{F}$ contains $n$ members in convex position.
 \end{lema}
It is easy to show Lemma 2 using Lemma 1 and Theorem 3.

\begin{teorema} \label{bound} $M(n) \le ({2n-5 \choose n-2}+1){2n-4 \choose n-2} +1.$
 \end{teorema}

\noindent \textbf{Proof.}  Consider a body $C_0$ intersecting the boundary of the convex hull of the family $\mathcal{F}$.  Let $l$ be a tangent line to $C_0$ such that all members of $\mathcal{F}$ lie completely on one side of $l$.  Then by rotating $l$ counterclockwise along the boundary of $C_0$, we order the members of $\mathcal{F}\setminus C_0$ in the order that $l$ meets each member.  Furthermore, we denote $l_i$ to be the tangent line of $C_0$ and $C_i \in F\setminus C_0$ when $l$ first meets $C_i$ in this rotation.  See Figure 2.

 \begin{figure}[h]
\includegraphics[width=220pt]{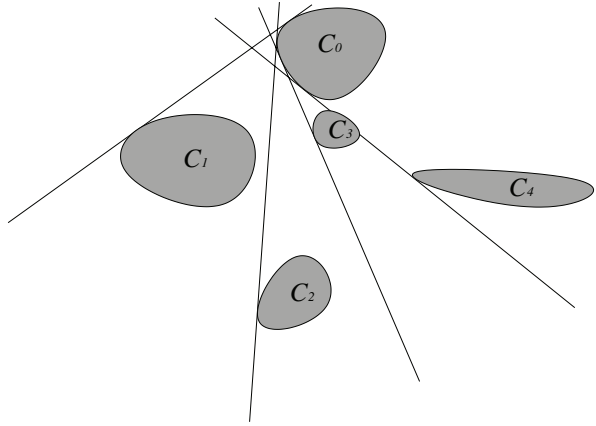}
 \caption{Ordering with tangent lines}
 \label{r}
\end{figure}

\noindent Set $N = {2n-5 \choose n-2}+1$.  By the pigeonhole principle, either

  \begin{enumerate}

  \item there exists a line that meets at least ${2n-4 \choose n-2} +1$ members of $\mathcal{F}$,

   \medskip

   \item or there exists convex bodies $ C_{i_1},C_{i_2},...,C_{i_N}$ and tangents lines $l_{i_1},l_{i_2},...,l_{i_N}$ such that the interior of $C_{i_k}$ does not intersect with any member of $\{l_{i_1},l_{i_2},...,l_{i_N}\}$ for all $1\leq k \leq N$.

       \end{enumerate}

 \begin{figure}
  \centering
    \subfigure[Case 1.]{\label{hide2}\includegraphics[width=.33\textwidth]{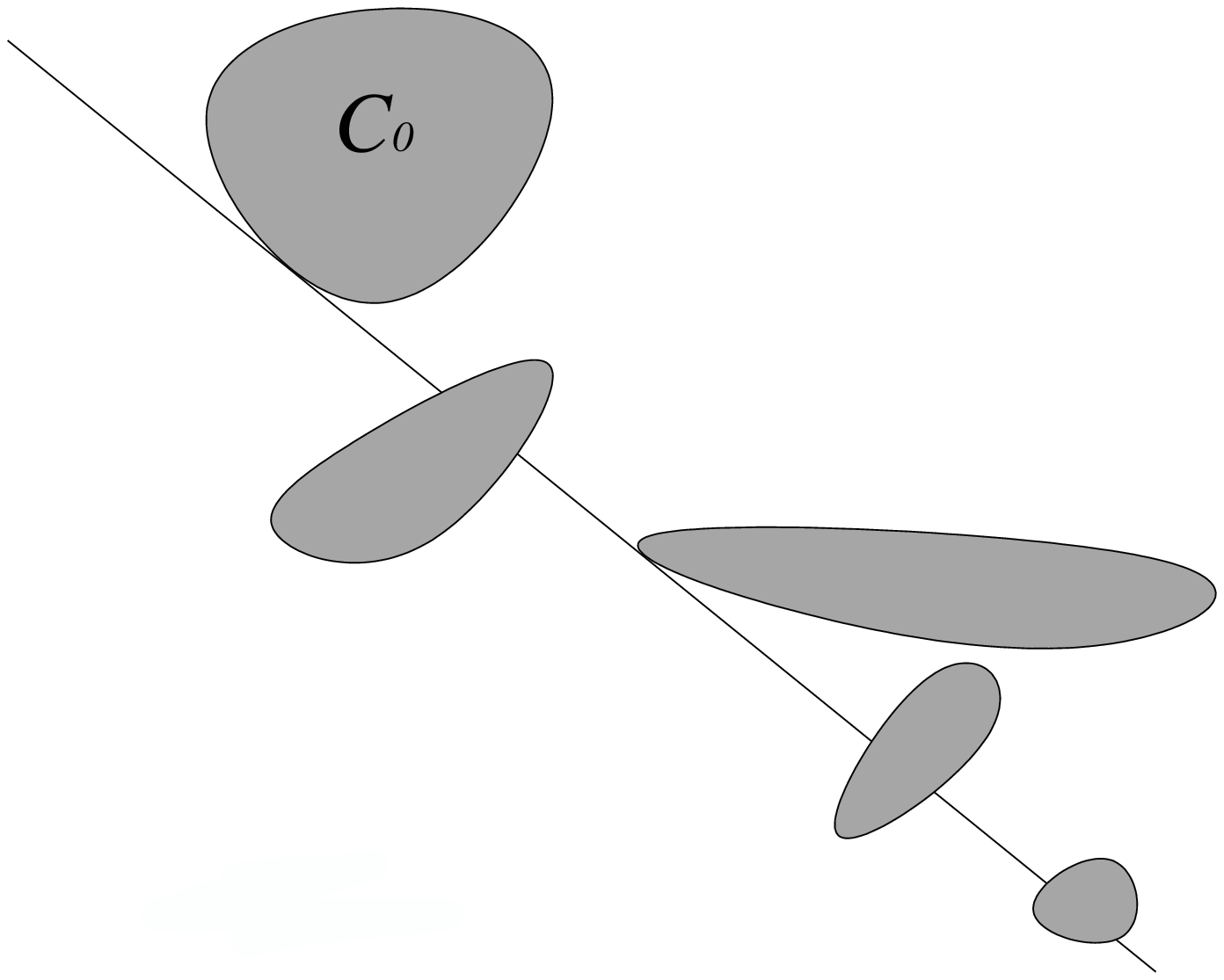}} \hspace{2cm}
      \subfigure[Case 2.]{\label{clockwise}\includegraphics[width=0.38\textwidth]{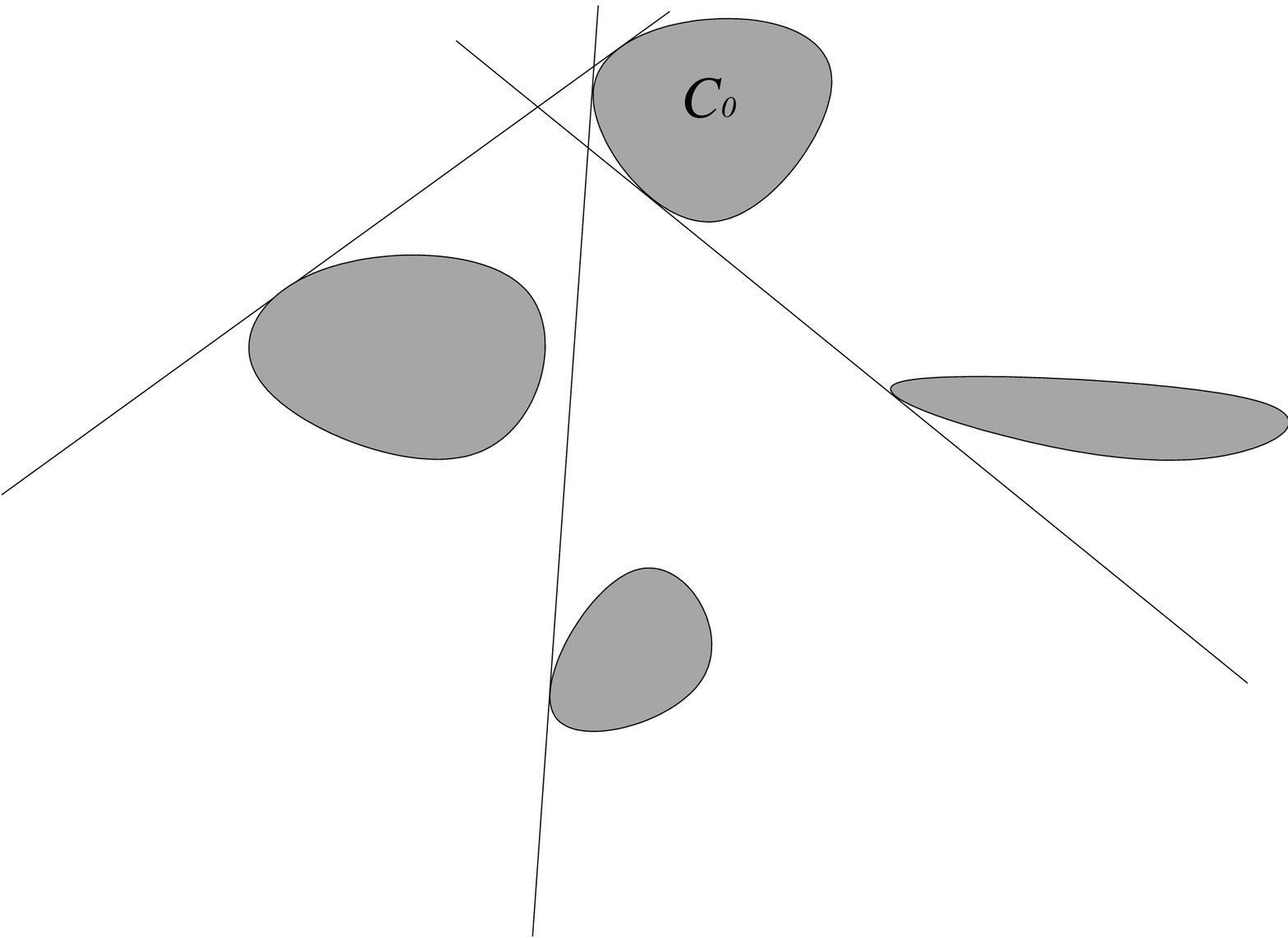}}

\caption{}
\label{fig:disconnect}
\end{figure}

\noindent See Figure 3.  In the first case, we can apply Lemma 2 to find $n$ members in convex position.  Therefore, suppose we are in the second case and let $\mathcal{F}_0  = \{C_{i_1},...,C_{i_N}\}$.

Assume that the triples $(C_{i_1},C_{i_2},C_{i_3})$ and $(C_{i_2},C_{i_3},C_{i_4})$ are oriented counterclockwise (clockwise).  We will show that $(C_{i_1},C_{i_2},C_{i_4})$ and $(C_{i_1},C_{i_3},C_{i_4})$ also have counterclockwise (clockwise) orientations.
We need to expose points $c_i \in C_i \cap bd[\conv[C_1,C_2,C_4]]$ such that $c_1,c_2,c_4$ are counterclockwise oriented.
Let $h_t$ denote a ray tangent to $C_0$, that starts at $C_0$ and goes as far as possible without leaving the convex hull of the family. The parameter $t$ increases as we move counterclockwise on the boundary of $C_0$. Let $c_1\in l_1 \cap C_1$. Let $c_2$ be contained in the intersection of $C_2$ with the exterior common tangent to $C_2$ and $C_3$ that has $conv[C_2,C_3]$ to it's right when oriented from $C_2$ to $C_3$. Finally let $c_4\in C_4 \cap h_{t^*_4}$, with $t^*_{4}$ the last $t$ such that $h_t \cap C_4\neq \emptyset$. Is easy to see that $\{c_1,c_2,c_4\}$ are counterclockwise oriented. Analogously $(C_1,C_3,C_4)$ is counterclockwise oriented. The clockwise case follows a similar argument that we skip.
Hence by Lemma 1, we can either find $n$ members in $\mathcal{F}_0$ such that ever triple has a clockwise orientation, or $n-1$ members such that every triple has a counterclockwise orientation.   Since $\{C_0,C_{i_j},C_{i_k}\}$ is counterclockwise oriented for every pair $j,k$, by Theorem 3 we can find $n$ members in convex position.\\

\begin{flushright}$ \blacksquare$\end{flushright}
\section{Order types}
Order types are natural combinatorial objects assigned to ordered point sets.  Given an ordered point set $X$ in $\mathbb{R}^d$ the order type can be defined as an orientation function $\chi: {X \choose d+1
} \to \{+,0,-\}$. See \cite{lectures}.

\subsection{On the Bisztriczky-Fejes T\'{o}th conjecture}
With theorem \ref{main} in mind we introduce the following definition.
\begin{definicion}
Let $\mathcal{F}$ be a family of noncrossing convex bodies in general position in the plane. Then $\mathcal{F}$ is said to be \textbf{representable} by the point set $X$ if there is a bijection $f:\mathcal{F} \to X$  such that, given any triple of points $Y \in {X \choose 3} $, if $Y$ is ordered so that  $\chi(Y)=+$, then $+ \in \chi(f^{-1}(Y))$, with $f^{-1}(Y)$ ordered by the pullback of the ordering of  $X$. \\
\end{definicion}

Remark: If a triple of convex bodies has both orientations then any triple of points represents the convex bodies.

\medskip

Take any representable family of $N(n)$ disjoint convex bodies in general position. Find a representation by points ($N(n)$ as in the Erd\H{o}s-Szekeres theorem). We can select a convex $n$-gon, and order this $n$-gon by orienting the boundary in the counterclockwise direction. The inverse image of this $n$-gon under $f$ with the induced (pulled back) ordering is in convex position by Theorem \ref{main}. So, if every disjoint family of convex bodies could be represented by points, the Bisztriczky-Fejes T\'{o}th conjecture would follow.\\

\begin{problem}[\bf{Hubard}] \label{rep}  Find the smallest integer $r(n)$ such that in every family of $r(n)$ disjoint convex bodies in general position there is a $n$-subfamily  that is representable by points.\end{problem}
The existence of $r(n)$ follows from the Bisztriczky-Fejes T\'{o}th theorem with $r(n) \le M(n)$.  J\'{a}nos Pach and Geza T\'{o}th  \cite{pach3}  a non-representable family of disjoint convex bodies. It consists of $9$ intervals and realizes a well known construction of a  non-stretchable pseudoline arrangement due to Ringel (sometimes denoted by Rin(9) or the non-Pappus configuration in the literature, see \cite{ormat}).  By flattening this construction to be contained in a neighborhood of an interval and iterating the construction, the authors \cite{pach3} were able to bound $r(n)$ by below. Their result yields the lower bound $r(n) > n^{\log9/\log8}$.\\
	
 A family of problems arise from this consideration. Can any family of convex bodies be represented by a family of intervals? More generally, what is the best representation of a family of convex bodies with property $\mathcal{P}$ by a family of convex bodies with property $\mathcal{Q}$? Here \emph{best representation} is defined as a Ramsey function similarly to $r(n)$. In \cite{suk} it was shown that if $\mathcal{P}$ is the class of families of segments in general position and $\mathcal{Q}$ is the class of point sets then $r_{\mathcal{P},\mathcal{Q}}(n+1) \leq n^4 +1$, i.e. any family of $n^4 +1$ disjoint segments in general position in the plane contains a subfamily of at least $n+1$ segments whose order type can be represented by points. \\

\subsection{Oriented Matroids}
As already mentioned, the example of \cite{pach3} is constructed realizing the order type of a non-stretchable pseudoline arrangement by a family of  segments.  A nice feature of the theory of oriented matroids is that many different structures turn out to be equivalent. For us the most natural approach is that of chirotopes (see \cite{ormat} for details).

We call a family of noncrossing convex bodies in general position \textbf{$3$-nondisconnectable} if for all $\mathcal{Y} \in {\mathcal{F} \choose 3}$, for all $A \in \mathcal{Y}$, $\conv(\mathcal{Y}\setminus A)$ is connected.\\ \\

\noindent\textbf{Remark} \emph{Every ordered $3$-nondisconnectable noncrossing family of convex bodies in general position $\mathcal{F}$ (with the function $\chi: {\mathcal{F} \choose 3} \to \{+,-\}$ defined as the orientation) forms a chirotope.}\\ \\
By corollary 3.6.3 in \cite{ormat} to the 3-term Grassman-Pl\"{u}cker relations, it suffices to show:

\begin{enumerate}

\item The mapping $\chi: {\mathcal{F} \choose 3} \to \{+,-\}$ is alternating and $\mathcal{F} \choose 3$ is the set of basis of a rank $3$ matroid.

\item The restriction of $\chi$ to any ${\mathcal{F} \choose 5}$ is realizable by points.

\end{enumerate}

\noindent Alternating means that for any triple and every $\sigma \in S_3$ (the symmetric group on three elements),
$\chi(A_1,A_2,A_3)= sign (\sigma)\chi(A_{\sigma(1)},A_{\sigma(2)},A_{\sigma(3)})$.\\
To prove 2) we only need to show that every family of $5$ convex bodies in the plane is realizable by points. This is easy to check by case analysis, splitting the cases by the number of convex bodies on the convex hull.

\begin{problem}[\bf{Hubard}]  Is every 3-uniform oriented matroid realizable by non-crossing convex bodies?
\end{problem}

One of the most beautiful problems in oriented matroid theory is the isotopy problem asked by Ringel in 1956:\\

\emph{ Given two point sets with the same order type is there a continuous path of point sets that goes from one to the other with the same order type at every moment?}\\ \\
This problem was solved in the negative by Mn\"ev and independently by several other researchers.   However Mn\"ev's Universality Theorem is the strongest result, see \cite{mnev}.

\begin{problem} [\bf{Hubard}]  Given two families with the same order type, is there a continuous path (under what topology?) of families that goes from one to the other with the same order type at every moment? \end{problem}

\section{Acknowledgments}
 The first three authors gratefully acknowledge the support of CONACYT and the SNI. The authors would also like to thank Andreas Holmsen, J\'anos Pach, Erik Dies, Javier Bracho and Imre B\'{a}r\'{a}ny. We will also like to dedicate this paper to the memory of Victor Neumann-Lara.\\
 \\
 \\
Alfredo Hubard\\
Courant Institute of Mathematical Sciences\\
 New York University\\
251 Mercer Street, New York, New York 10012\\
hubard@cims.nyu.edu\\
\\
Luis Montejano\\
Instituto de Matematicas\\
UNAM\\
Mexico DF 04510\\
montejano@matem.unam.mex\\
\\
Emiliano Mora\\
Instituto de Matematicas\\
UNAM\\
Mexico DF 04510\\
emailiano@gmail.com\\
 \\
Andrew Suk\\
Courant Institute, New York, New York and\\
EPFL, Lausanne, Switzerland\\
suk@cims.nyu.edu\\
\\

\end{document}